\documentclass[onefignum,onetabnum]{siamart190516}

\usepackage{lipsum}
\usepackage{amsfonts}
\usepackage{graphicx}
\usepackage{epstopdf}
\usepackage{algorithmic}
\ifpdf
  \DeclareGraphicsExtensions{.eps,.pdf,.png,.jpg}
\else
  \DeclareGraphicsExtensions{.eps}
\fi

\newsiamremark{hypothesis}{Hypothesis}
\crefname{hypothesis}{Hypothesis}{Hypotheses}
\newsiamthm{claim}{Claim}

\headers{Adaptive Ritz Method
}{M. Liu and Z. Cai}

\title{Adaptive Two-Layer ReLU Neural Network:\\
II. Ritz Approximation to Elliptic PDEs\thanks{This work was supported in part by the National Science Foundation
under grant DMS-2110571.}}

\author{Min Liu\thanks{School of Mechanical Engineering, Purdue University, 585 Purdue Mall, West Lafayette, IN 47907-2088(\email{liu66@purdue.edu}). }
\and Zhiqiang Cai\thanks{Department of Mathematics, Purdue University, 150 N. University Street, West Lafayette, IN 47907-2067 
  (\email{caiz@purdue.edu}).}
}

\usepackage{amsopn}

\ifpdf
\hypersetup{
  pdftitle={FOSLS},
 pdfauthor={M.Liu, Z. Cai, and R.Falgout}
}
\fi

\usepackage[utf8]{inputenc}
\usepackage{bm}
\usepackage{mathtools}
\usepackage{indentfirst}
\usepackage{subfigure}
\usepackage{tcolorbox}
\usepackage{color}
\usepackage[export]{adjustbox}
\usepackage{scalerel}

\usepackage{algorithmic}
\Crefname{ALC@unique}{Line}{Lines}

\newcommand{\R}{\mathbb{R}}

\usepackage{graphicx}
\usepackage{diagbox}
\usepackage{pifont}
\newcommand{\vertiii}[1]{{\left\vert\kern-0.25ex\left\vert\kern-0.25ex\left\vert #1 
    \right\vert\kern-0.25ex\right\vert\kern-0.25ex\right\vert}}

\newcommand{\bsigma}{\mbox{\boldmath${\sigma}$}}

\newcommand{\btheta}{\mbox{\boldmath${\theta}$}}

\newcommand{\bomega}{\mbox{\boldmath${\omega}$}}
\newcommand{\btau}{\mbox{\boldmath${\tau}$}}

\newcommand{\bvarphi}{\mbox{\boldmath${\varphi}$}}

\newcommand{\bn}{{\bf n}}

\newcommand{\bc}{{\bf c}}
\newcommand{\be}{{\bf e}}
\newcommand{\bl}{{\bf l}}

\newcommand{\bx}{{\bf x}}
\newcommand{\bv}{{\bf v}}

\newcommand{\cK}{{\cal K}}
\newcommand{\cM}{{\cal M}}

\newcommand{\cS}{{\cal S}}
\newcommand{\cT}{{\cal T}}

\begin{document}

\maketitle
\begin{abstract}
In this paper, we study adaptive neuron enhancement (ANE) method for solving self-adjoint second-order elliptic partial differential equations (PDEs). The ANE method is a self-adaptive method generating a two-layer spline NN and a numerical integration mesh such that the approximation accuracy is within the prescribed tolerance. Moreover, the ANE method provides a natural process for obtaining a good initialization which is crucial for training nonlinear optimization problem.

The underlying PDE is discretized by the Ritz method using a two-layer spline neural network based on either the primal or dual formulations that minimize the respective energy or complimentary functionals. Essential boundary conditions are imposed weakly through the functionals with proper norms. It is proved that the Ritz approximation is the best approximation in the energy norm; moreover, effect of numerical integration for the Ritz approximation is analyzed as well. Two estimators for adaptive neuron enhancement method are introduced, one is the so-called recovery estimator and the other is the least-squares estimator. Finally, numerical results for diffusion problems with either corner or intersecting interface singularities are presented. 
\end{abstract}



\section{Introduction}

Recent success of neural networks (NNs) for many artificial intelligence tasks has led wide applications to other fields, including recent studies of using neural network models to numerically solve partial differential equations (PDEs) (see, e.g., \cite{Berg18, CAI2020, Weinan18, Karniadakis19, Sirignano18}). 
Because neural network functions are nonlinear functions of the parameters, 
discretization of a PDE is set up as an optimization problem through either the natural minimization or manufactured least-squares (LS) principles. Hence, existing methods consist of (1) the deep Ritz method \cite{Weinan18} and (2) the deep LS method such as the deep Galerkin method (DGM) \cite{Sirignano18}, the physics-informed neural networks (PINN) \cite{Karniadakis19}, the deep LS and FOSLS methods \cite{CAI2020}, etc.  
The former has the least variables, requires the least smoothness, but is applicable to a single class of problems. The latter is applicable to a large class of PDEs, but either has additional variables such as the FOSLS or requires additional smoothness like the LS. 

Neural networks produce a new class of functions through compositions of linear transformations and activation functions. This class of functions is extremely rich. For example, it contains piecewise polynomials, which are the footing of spectral elements, and continuous and discontinuous finite element methods. It approximates polynomials of any degree with exponential efficiency, even using simple activation functions like ReLU.
Despite many efforts, it is widely accepted that approximation properties of NNs are not yet well-understood. 
As a consequence, design of network structures for approximating the solution of a PDE within the prescribed accuracy remains open and is mainly done by time consuming trial-and-error. 


To address the issue on how to design a minimal network model required to approximate a function/PDE within the prescribed accuracy, in terms of width, depth, and the number of parameters, we propose and study the adaptive neuron enhancement (ANE) method for approximating a function in \cite{LiuCai} and for solving a self-adjoint second-order elliptic PDE in this paper.  Specifically, for a given tolerance $\epsilon>0$, the ANE method generates a two-layer spline neural network such that the approximation accuracy is within the prescribed tolerance. 
The key ingredient of the method is the neuron enhancement strategy which determines how many new neurons to be added, when the current approximation is not within the given accuracy. This is done through local error indicators collected on the physical subdomains, and a proper neuron initialization.

The underlying PDE is discretized by the Ritz method using a two-layer spline neural network based on either the primal or dual formulations. The primal problem minimizes the energy functional with the (essential) Dirichlet boundary condition imposed weakly \cite{CAI2020}. Another way to impose the Dirichlet boundary condition is the well-known Nitsche’s method that requires the stabilization constant is sufficient large (see \cite{Xu2020} in the context of neural networks).
When the dual (flux) variable is important for the underlying application, we may maximize the complementary functional directly; in this case, the Neumann boundary condition becomes essential and is again enforced weakly through the functional with a proper norm. 

Even though the set of neural network functions does not form a space, we show that the Ritz approximation based on either the primal or the dual formulation is the best approximation in the energy norm. Moreover, we are able to analyze the effect of numerical integration as well. As expected, the error in the energy norm by the Ritz approximation with numerical integration is bounded by the approximation error of the neural network and the error of the numerical integration. This result may be considered as the extension of the first Strang's lemma (see, e.g., \cite{Ciarlet78}) from the Galerkin approximation over a subspace to the Ritz approximation over a set.

{\it A posteriori} error estimation is important for the ANE method. 
In this paper, we consider two estimators: the recovery and the least-squares estimators. For the primal formulation, the recovery estimator is defined as a weighted $L^2$ norm of difference between the numerical and the recovered fluxes. When the recovered flux is more accurate than the numerical flux (i.e., the so-called saturation assumption \cite{Bank}), the recovery estimator is proved to be reliable and efficient. The least-squares estimator \cite{cai2010} adds an additional term, the $L^2$ norm of the residual, to the recovery estimator.

The paper is organized as follows.  
The primal and dual formulations of the diffusion-reaction problem and their well-posedness are discussed in section 2. The Ritz approximation using neural networks are described in section 3. Error estimates of the Ritz approximations and effect of numerical integration are obtained in sections 3 and 4, respectively. {\it A posteriori} error estimators, the ANE method, and initialization at each stage of the ANE method are introduced in the respective section 5, 6, and 7. Finally, we present numerical results for problems with either corner or intersecting interface singularities in section 8. 



In this paper, we will use the standard notation and definitions for the Sobolev space $H^{s}(\Omega)$ and $H^{s}(\Gamma)$ for a subset $\Gamma$ in $\partial \Omega$. 
The standard associated inner product and norms are denoted by $(\cdot,\cdot)_{s,\Omega}$ and $(\cdot,\cdot)_{s,\Gamma}$ and by $\Vert \cdot \Vert _{s,\Omega}$ and $\Vert \cdot \Vert _{s,\Gamma} $, respectively. 
When $s=0$, $H^{0}(\Omega)$ coincides with $L^2(\Omega)$.
Denote the corresponding norms on product space $H^s(\Omega)^d$ by $\|\cdot\|_{s,\,\Omega,\,d}$ and
$|\cdot|_{s,\,\Omega,\,d}$.  When there is no ambiguity, the subscript $\Omega$ and $d$ in the designation of norms  will be suppressed. 

\section{Diffusion-Reaction Problem}

Let $\Omega$ be a bounded domain in ${\R}^d$ with
Lipschitz boundary $\partial\Omega= \bar{\Gamma}_D \cup \bar{\Gamma}_N$. 
Consider the following self-adjoint second-order scalar elliptic partial
differential equation:
 \begin{equation}\label{pde}
 \left\{\begin{array}{rlll}
 -\mbox{div}\, (A \nabla \,u) +c\,u & = & f, & \mbox{in }  \Omega\subset \R^d, \\[2mm]
 u &=&  g_{_{\small D}}, &  \mbox{on} \ \Gamma_D, \\[2mm]
 -{\bf n} \cdot A \nabla \, u &=& g_{_{\small N}}, & \mbox{on} \ \Gamma_N,
 \end{array}\right.
\end{equation}
where $f \in L^2(\Omega)$, $c \in L^2(\Omega)$, $g_{_{\small D}}\in H^{1/2}(\Gamma_D)$, $g_{_{\small N}}\in H^{-1/2}(\Gamma_N)$;
$A(x)$ is a $d \times d$ symmetric
matrix-valued function in $L^2(\Omega)^{d\times d}$; 
and ${{\bf n}}$ is the outward
unit vector normal to the boundary. We assume that $A$ is uniformly positive definite and that $c\ge 0$. 


The natural optimization problem of (\ref{pde}) is the so-called primal problem that minimizes the energy functional over $H^1_{g,D}(\Omega)=\{v\in H^1(\Omega):\, v=g_D\mbox{ on } D\}$.
Since it is difficult for neural network functions to satisfy boundary conditions (see \cite{Weinan18}), as in \cite{CAI2020}, we enforce the Dirichlet (essential) boundary condition weakly through the energy functional. 
To this end, define the energy functional by
\begin{eqnarray}\nonumber
J(v) &=&  \dfrac12\left\{ \|A^{1/2}\nabla v\|^2_{0,\Omega}
+\|c^{1/2}v\|^2_{0,\Omega}+ \gamma_D \|v-g_D\|^2_{1/2,\Gamma_D}\right\} -\Big((f,v) +(g_N, v)_{0,\Gamma_N}\Big)\\[2mm] \label{energy}
&=& \dfrac12 a(v,v) -f(v) + \dfrac{\gamma_D}{2} \|g_D\|^2_{1/2,\Gamma_D},
\end{eqnarray}
where $\gamma_D>0$ is a constant and the quadratic form $a(\cdot,\cdot)$ and the linear form $f(\cdot)$ are given by
 \[
 a(v,v)= \|A^{1/2}\nabla v\|^2_{0,\Omega}
+\|c^{1/2}v\|^2_{0,\Omega}+ \gamma_D \|v\|^2_{1/2,\Gamma_D}
\]
and
\[
f(v) = (f,v) +(g_N, v)_{0,\Gamma_N} + \gamma_D (g_D,v)_{1/2,\Gamma_D},
 \]
respectively. 
Then the primal problem is to find $u\in V:= H^1(\Omega)$ such that
 \begin{equation}\label{primal}
 J(u) = \min_{v\in V} J(v).
 \end{equation}
 
\begin{remark}
Another way to enforce the Dirichlet boundary condition is the well-known Nitsche’s method usually stated in the Galerkin formulation. Its equivalent form for the Ritz formulation is to transform the Dirichlet boundary condition to the Robin boundary condition by penalization {\em (}see {\em \cite{Xu2020}} in the context of neural networks{\em )}, and the penalization constant is usually required to be large.
\end{remark}

\begin{proposition}
Problem {\em (\ref{primal})} has a unique solution $u\in V$. Moreover, the solution $u$ satisfies the following {\it a priori} estimate:
\begin{equation}\label{stable}
\|u\|_{1,\Omega}\leq C\, \left( \|f\|_{-1,\Omega} + \|g_D\|_{1/2,\Gamma_D}
+ \|g_N\|_{-1/2,\Gamma_N}\right).
\end{equation}
\end{proposition}

\begin{proof}
By the assumptions on $A$ and $c$ and the trace theorem, the bilinear form $a(\cdot,\cdot)$ is $H^1(\Omega)$ elliptic, i.e., there exists a positive constant $\alpha$ such that 
 \begin{equation}\label{ellipticity}
 \alpha \|v\|^2_{1,\Omega} \leq a(v,v) , \quad\forall\,\, v\in V.
 \end{equation}
It is easy to see that the bilinear form $a(\cdot,\cdot)$ and the linear form $f(\cdot)$ are continuous in $V\times V$ and $V$, respectively. Then the Lax-Milgram lemma implies that problem (\ref{primal}) has one and only one solution in $V$. 

The solution $u\in V$ of problem (\ref{primal}) may be characterized by the relation
 \[
 a(u,v)=f(v), \quad\forall\,\, v\in V.
 \]
Now, the {\it a priori} estimate in (\ref{stable}) is a direct consequence of (\ref{ellipticity}) and the fact that 
\[
|f(u)|\leq C \,  \left( \|f\|_{-1,\Omega} + \|g_D\|_{1/2,\Gamma_D}
+ \|g_N\|_{-1/2,\Gamma_N}\right) \, \|u\|_{1,\Omega}.
\]
This completes the proof of the proposition. 
\end{proof}

Another optimization problem of (\ref{pde}) is the so-called dual problem that maximizes the complementary functional for the dual variable $\bsigma = - A\nabla u$ over the dual space
\[
\Sigma:=H(\mbox{div};\Omega) =\{\btau\in L^2(\Omega)^d : \, \nabla\cdot \btau \in L^2(\Omega)\}.
\]
For simplicity of presentation, assume that the diffusion coefficient $c(\bx)$ is positive.
The negative of the complementary functional is given by 
\begin{eqnarray}\nonumber 
    J^*(\btau) &=&  \dfrac12\left\{ \|A^{-1/2}\btau\|^2_{0,\Omega}+\|\dfrac{1}{\sqrt{c}}(\nabla\!\cdot\!\btau -f)\|^2_{0,\Omega} +\gamma_{\small N}\|\btau\!\cdot\!\bn+g_N\|^2_{-1/2,\Gamma_N}\! \right\}  +\! \int_{\Gamma_D}\!\! g_D (\btau\!\cdot\!\bn)ds
    \\[2mm] \label{complementary}
    &=& \dfrac12 a^*(\btau,\btau) -f^*(\btau) +\dfrac12 \big(\|c^{-1/2}f\|^2_{0,\Omega} + \gamma_N \|g_N\|^2_{-1/2,\Gamma_N}
    \big),
\end{eqnarray}
where $\gamma_N>0$ is a constant and the quadratic form $a^*(\cdot,\cdot)$ and the linear form $f^*(\cdot)$ are given by
 \[
 a(\btau,\btau)= \|A^{-1/2}\btau\|^2_{0,\Omega}+\|c^{-1/2}\nabla\!\cdot\!\btau \|^2_{0,\Omega} +\gamma_{\small N}\|\btau\!\cdot\!\bn\|^2_{-1/2,\Gamma_N}\!
\]
and
\[
f^*(\btau) = (f,c^{-1}\nabla\!\cdot\!\btau) -(g_D, \btau\!\cdot\!\bn)_{0,\Gamma_D} - \gamma_N (g_N,\btau\!\cdot\!\bn)_{-1/2,\Gamma_N},
 \]
respectively. 
The Neumann boundary condition becomes essential for the dual problem and is enforced weakly through the complementary functional defined in (\ref{complementary}).
The dual problem is then to seek $\bsigma\in \Sigma$ such that
 \begin{equation}\label{dual}
 J^*(\bsigma) = \min_{\btau\in \Sigma} J^*(\btau). 
 \end{equation}
 
\begin{proposition} Problem {\em (\ref{dual})} has a unique solution $\bsigma\in \Sigma$. Moreover, the solution $\bsigma$ satisfies the following {\it a priori} estimate:
\[\|\bsigma\|_{0,\Omega}\leq C\, \left( \|f\|_{-1,\Omega} + \|g_D\|_{1/2,\Gamma_D}+ \|g_N\|_{-1/2,\Gamma_N}\right).\]
\end{proposition}

\begin{proof}
The proposition may be proved in a similar fashion as that of Proposition~2.1.
\end{proof}

\section{Neural Network Methods}

A two-layer neural network (NN) consists of an input and output layers. The output layer usually has no activation function. Let $\tau_k$ be the spline activation function of the form:
 \[
 \tau_k(t) = \max \{0,\,t^k\}
 =\left\{\begin{array}{ll}
   0, &  t < 0 ,\\[2mm]
   t^k , & t \ge 0,
   \end{array}\right.
 \]
for a fixed integer $k>0$. It is clear that $\tau_k(t)$ is a piece-wise polynomial of degree $k$ in $C^{k-1}(\R)$ with one breaking point $t=0$. When $k=1$, the activation function $\tau_1(t)$ is the popular rectified linear unit (ReLU).
A two-layer spline NN with $n$ neurons generates the following functional class from $\R^d\rightarrow \R^o$:
\begin{equation}\label{ReLU-n}
{\cal M}^o_n(\tau_k) = \left\{\bc_0+\sum_{i=1}^n \bc_i\tau_k(\bomega_i\!\cdot\! \bx -b_i)\, :\,  
 b_i\in \R,\,\, \bomega_i\in \cS^{d-1},\,\, \bc_i\in \R^o \right\},
 \end{equation}
where $d$ and $o$ are the respective dimension of the input and output; $\bomega_i\in \cS^{d-1}$ 
and $b_i\in \R$ are the respective weights and bias of the input layer; $\bc_i\in\R^o$ and $\bc_0\in\R^o$ are the respective weights and bias of the output layer; and $\cS^{d-1}$ is the unit sphere in $\R^d$. The total number of parameters of ${\cal M}^o_n(\tau_k)$ is
\[
M_d(n)=(d+o)n+o,
\]
where $\{\bc_i\}_{i=0}^n$ are linear parameters and $\{\bomega_i, b_i\}_{i=1}^n$ are nonlinear parameters.

\begin{remark}
Since $\tau_k$ is not a polynomial,
it has been proven (see, e.g., \cite{Cybenko1989, HornikS1989}) that the linear space
\[
{\cal M}(\tau_k) = 
\left\{v(\bx)\in {\cal M}^1_n(\tau_k):\,  
n\in \mathbb{Z}_+ \right\}
\]
is dense in $C(K)$, the space of all continuous functions defined on a compact set $K\in \R^d$.
Note that ${\cal M}^o_n(\tau_k)$ is a subset, but not a subspace, of ${\cal M}^o(\tau_k)=\underbrace{{\cal M}(\tau_k)\times \cdots \times {\cal M}(\tau_k)}_o$. 

Even though results on approximation order are still scarce, there are two noticeable results for target functions in Sobolev space in the $L^2(\Omega)$ norm by Petrushev {\em \cite{Petrushev1998}} and in spectral Barron space in the Sobolev norm by Siegel and Xu {\em \cite{siegel2020approximation2}}. It is not clear if the former has been extended to the $H^1(\Omega)$ norm. For the latter, solutions of very few partial differential equations have been shown in the spectral Barron space {\em \cite{Weinan20Barron}}.
\end{remark}

To approximate the solution of (\ref{pde}) using neural network functions, the Ritz method is to minimize the energy functional over the set $\cM^1_n(\tau_k)$, i.e., 
finding $u_n\in \cM^1_n(\tau_k)\subset H^1(\Omega)$ such that 
\begin{equation}\label{Ritz}
J(u_n) = \min_{v\in \cM^1_n(\tau_k)} J(v).
\end{equation}
Since $\cM^1_n(\tau_k)$ is not convex, problem (\ref{Ritz}) has many solutions.  

\begin{theorem}
Let $u\in H^1(\Omega)$ be the solution of problem {\em (\ref{primal})}, and let $u_n\in \cM^1_n(\tau_k)$ be a solution of {\em (\ref{Ritz})}. Then we have
 \begin{equation}\label{error-R}
     \|u-u_n\|_a = \inf_{v\in \cM^1_n(\tau_k)} \|u-v\|_a,
 \end{equation}
where $\|v\|_a:=\sqrt{a(v,v)}$ is the energy norm for the primal variable. 
\end{theorem} 

\begin{proof}
The proof of the theorem follows easily from the standard error estimate for the Rite approximation: for any $v\in \cM^1_n(\tau_k)$,
\[
\|u-u_n\|_a^2=2\left(J(u_n)-J(u)\right)
\leq 2\left(J(v)-J(u)\right) = \|u-v\|_a^2,
\]
which implies the validity of (\ref{error-R}). 
\end{proof}

When the flux variable $\bsigma=-A\nabla u$ is important for the underlying application, we may approximate it directly through the dual problem: finding
$\bsigma_n\in \cM^d_n(\tau_k)$ such that 
\begin{equation}\label{Dual-n}
J^*(\bsigma_n) = \min_{\btau\in \cM^d_n(\tau_k)} J^*(\btau).
\end{equation}

\begin{theorem}
Let $\bsigma\in H(\mbox{div};\Omega)$ be the solution of problem {\em (\ref{dual})}, and let $\bsigma_n\in \cM^d_n(\tau_k)$ be a solution of {\em (\ref{Dual-n})}. Then we have
 \begin{equation}\label{error-D}
     \|\bsigma-\bsigma_n\|_{a^*} = \inf_{\btau\in \cM^d_n(\tau_k)} \|\bsigma-\btau\|_{a^*},
 \end{equation}
where $\|\btau\|_{a^*}:=\sqrt{a^*(\btau,\btau)}$ is the energy norm for the dual variable. 
\end{theorem} 

\begin{proof}
The theorem may be proved in a similar fashion. For any $v\in \cM^1_n(\tau_k)$,
\[
\|\bsigma-\bsigma_n\|_{a^*}^2=2\left(J^*(\bsigma_n)-J^*(\bsigma)\right)
\leq 2\left(J^*(\btau)-J^*(\bsigma)\right) = \|\bsigma-\btau\|_{a^*},
\]
which implies the validity of (\ref{error-D}).
\end{proof}

\section{Effect of Numerical Integration}

In practice, the integrals of the loss functional are computed numerically. For example, we proposed and implemented the composite mid-point quadrature rule in \cite{CAI2020}. To understand the effect of numerical integration, we extend the first Strang lemma for the Galerkin approximation over a subspace (see, e.g, \cite{Ciarlet78}) to the Ritz approximation over a subset.

To this end, denote by $a_{_\cT}(\cdot,\cdot)$ and $f_{_\cT}(\cdot)$ the discrete counterparts of $a(\cdot,\cdot)$ and $f(\cdot)$ through numerical integration. Similarly, $a^*_{_\cT}(\cdot,\cdot)$ and $f^*_{_\cT}(\cdot)$ are the discrete counterparts of $a^*(\cdot,\cdot)$ and $f^*(\cdot)$. Then approximations to the primal and dual variables with numerical integration are seeking $u_{_\cT}\in \cM^1_n(\tau_k)$ such that 
\begin{equation}\label{Ritz-n-Q}
J_{_\cT}(u_{_\cT}) = \min_{v\in \cM^1_n(\tau_k)} J_{_\cT}(v), \quad\mbox{where }\,\, J_{_\cT}(v)=\dfrac12 a_{_\cT}(v,v) -f_{_\cT}(v)
\end{equation}
and $\bsigma_{_\cT}\in \cM^d_n(\tau_k)$ such that 
\begin{equation}\label{Dual-n-Q}
J^*_{_\cT}(\bsigma_{_\cT}) = \min_{\btau\in \cM^d_n(\tau_k)} J^*_{_\cT}(\btau), \quad\mbox{where }\,\, J^*_{_\cT}(v)=\dfrac12 a_{_\cT}(\btau,\btau) -f_{_\cT}(\btau),
\end{equation}
respectively.


\begin{theorem}
Assume that there exists a positive constant $\alpha$ independent of $\cM^1_{2n}(\tau_k)$ such that
 \begin{equation}\label{Q}
     \alpha\, \|v\|_a^2 \leq a_{_\cT}(v,v), \quad\forall\,\, v\in \cM^1_{2n}(\tau_k).
 \end{equation}
Let $u$ and $u_{_\cT}$ be the solutions of {\em (\ref{primal})} and {\em (\ref{Ritz-n-Q})}, respectively. Then there exists a positive constant $C$ such that 
\begin{eqnarray}\nonumber   
    && \quad  \|u-u_{_\cT}\|_a
     \\[2mm] \label{E-Q}
    \qquad\quad &\leq & C\left(\inf_{v\in \cM^1_n(\tau_k)} \left\{\|u-v\|_a + \sup _{\phi\in \cM^1_{2n}(\tau_k)} \dfrac{|a(v,\phi)-a_{_\cT}(v,\phi)|}{\|\phi\|_a}\right\}+ \sup _{\phi\in \cM^1_{2n}(\tau_k)} \dfrac{|f(\phi) - f_{_\cT}(\phi)|}{\|\phi\|_a} 
    \right).
\end{eqnarray}
\end{theorem}

\begin{proof}
For any $v\in \cM^1_{n}(\tau_k)$, it is easy to see that $u_{_\cT}-v\in \cM^1_{2n}(\tau_k)\subset V$.
By the assumption in (\ref{Q}), the definition of $J_{_\cT}(\cdot)$, and the relations:
 \[
 J_{_\cT}(u_{_\cT}) \leq J_{_\cT}(v)
 \quad\mbox{and}\quad
 a(u,u_{_\cT}-v)=f(u_{_\cT}-v),
 \] 
we have
\begin{eqnarray*}
\dfrac{\alpha}{2} \|u_{_\cT}-v\|^2_a
 &\leq &\dfrac12 a_{_\cT}(u_{_\cT}-v,u_{_\cT}-v)
 = J_{_\cT}(u_{_\cT}) -J_{_\cT}(v) +f_{_\cT}(u_{_\cT}-v) -a_{_\cT}(v,u_{_\cT}-v)  \\[2mm]
&\leq & f_{_\cT}(u_{_\cT}-v)-a_{_\cT}(v,u_{_\cT}-v)  \\[2mm] \nonumber
&=& \Big(f_{_\cT}(u_{_\cT}-v) -f(u_{_\cT}-v)\Big) +\Big(a(v,u_{_\cT}-v)-a_{_\cT}(v,u_{_\cT}-v)\Big) + a(u-v, u_{_\cT}-v)
\end{eqnarray*}
which, together with the Cauchy-Schwarz inequality, implies
\[
    \|u_{_\cT}-v\|^2_a
    \leq C \left(\|u-v\|^2_a + \sup _{\phi\in \cM^1_{2n}(\tau_k)} \dfrac{|a(v,\phi)-a_{_\cT}(v,\phi)|}{\|\phi\|_a}+ \sup _{\phi\in \cM^1_{2n}(\tau_k)} \dfrac{|f(\phi) - f_{_\cT}(\phi)|}{\|\phi\|_a} 
    \right).
\]
Combining the above inequality with the triangle inequality
 \[
 \|u-u_{_\cT}\|_a \leq \|u-v\|_a + \|v-u_{_\cT}\|_a
 \]
and taking the infimum over all $v\in \cM^1_{n}(\tau_k)$ yield (\ref{E-Q}). This completes the proof of the theorem.
\end{proof}

\begin{remark}
For any $\phi_1,\, \phi_2\in \cM^1_{n}(\tau_k)$, since $\cM^1_{n}(\tau_k)$ is not a subspace, $\phi_1-\phi_2$ is generally not in $\cM^1_{n}(\tau_k)$. But it is easy to see that $\phi_1-\phi_2\in \cM^1_{2n}(\tau_k)$. This is why
the assumption in {\em (\ref{Q})} and the supremum in {\em (\ref{E-Q})} are over $\cM^1_{2n}(\tau_k)$ but not $\cM^1_{n}(\tau_k)$. 
\end{remark}

\begin{theorem}
Assume that there exists a positive constant $\alpha^*$ independent of $\cM^d_{2n}(\tau_k)$ such that
 \begin{equation}\label{Q1}
     \alpha^*\, \|\btau\|_{a^*}^2 \leq a^*_{_\cT}(\btau,\btau), \quad\forall\,\, \btau\in \cM^d_{2n}(\tau_k).
 \end{equation}
Let $\bsigma$ and $\bsigma_{_\cT}$ be the solutions of {\em (\ref{dual})} and {\em (\ref{Dual-n-Q})}, respectively. Then there exists a positive constant $C$ such that 
\begin{eqnarray}\nonumber   
    && \quad  \|\bsigma-\bsigma_{_\cT}\|_{a^*}
     \\[2mm] \label{E-Q-sigma}
    \quad\quad &\leq & C\left(\!\inf_{\btau\in \cM^d_n(\tau_k)}\!\! \left\{\|\bsigma-\btau\|_{a^*} + \!\! \sup _{\bv\in \cM^d_{2n}(\tau_k)} \!\! \dfrac{|a(\btau,\bv)-a_{_\cT}(\btau,\bv)|}{\|\bv\|_{a^*}}\right\}+ \!\!\sup _{\bv\in \cM^d_{2n}(\tau_k)} \!\! \dfrac{|f(\bv) - f_{_\cT}(\bv)|}{\|\bv\|_{a^*}} 
    \right).
\end{eqnarray}
\end{theorem}

\begin{proof}
The theorem may be proved in a similar fashion as that of Theorem 4.2.
\end{proof}

\section{Estimators}

Like adaptive finite element method, {\it A posteriori} error estimation plays a crucial role for the ANE method. The {\it a posteriori} error estimator is used 
for determining if the current approximation is within the target accuracy and {\it a posteriori} error indicators for determining how many new neurons to be added, where to refine the integration mesh, and how to initialize parameters of new neurons.
By following ideas of the {\it a posteriori} error estimation for the finite element method (see, e.g., \cite{Verfurth2013}),
we study the recovery and least-squares estimators in this section. 

Let $u_{_\cT}\in \cM^1_n(\tau_k)$ be a solution of (\ref{Ritz-n-Q}). Define the recovered flux $\hat{\bsigma}_{_\cT}\in \cM^d_n(\tau_k)$ satisfying
 \begin{equation}\label{flux}
     \|D^{-1/2}\big(\hat{\bsigma}_{_\cT} + A\nabla u_{_\cT}\big)\|_{0,\Omega} = \inf_{\btau\in \cM^d_n(\tau_k)} \|D^{-1/2}\big(\btau + A\nabla u_{_\cT}\big)\|_{0,\Omega} ,
 \end{equation}
where $D$ is either the identity $I$, $A$, or $A^2$.
Then the estimator and indicators are given by 
 \begin{equation}\label{xi}
     \xi = \|D^{-1/2}\big(\hat{\bsigma}_{_\cT} + A\nabla u_{_\cT}\big)\|_{0,\Omega}
     \quad\mbox{and}\quad
     \xi_K = \|D^{-1/2}\big(\hat{\bsigma}_{_\cT} + A\nabla u_{_\cT}\big)\|_{0,K}, \quad\forall\,\, K\in \cK_n,
 \end{equation}
where $\cK_n = \{K\}$ is the physical partition of the domain $\Omega$ for the current approximation $u_{_\cT}$ (see the subsequent section). To analyze the estimator $\xi$, we make the standard saturation assumption (see, e.g., \cite{Bank}): there exists a positive constant $\gamma\in [0,1)$ such that
 \begin{equation}\label{satu}
     \|D^{-1/2}\big(\hat{\bsigma}_{_\cT} + A\nabla u\big)\|_{0,\Omega}
 \leq \gamma\, \|D^{-1/2}A\nabla \big(u-u_{_\cT}\big)\|_{0,\Omega}.
 \end{equation}

\begin{theorem}
Under the assumption in {\em (\ref{satu})}, we have
\begin{equation}\label{xi-b}
    \dfrac{1}{1+\gamma}\, \xi \leq  \|D^{-1/2}A\nabla \big(u-u_{_\cT}\big)\|_{0,\Omega} \leq \dfrac{1}{1-\gamma}\, \xi.
\end{equation}
\end{theorem} 

\begin{proof}
The first inequality in (\ref{xi-b}) is a direct consequence of the triangle inequality and (\ref{satu}):
 \[
 \xi \leq \|D^{-1/2}\big(\hat{\bsigma}_{_\cT} + A\nabla u\big)\|_{0,\Omega} + \|D^{-1/2}A\nabla \big(u-u_{_\cT}\big)\|_{0,\Omega}
 \leq (1 +\gamma) \, \|D^{-1/2}A\nabla \big(u-u_{_\cT}\big)\|_{0,\Omega}.
 \]
The second inequality in (\ref{xi-b}) may be proved in a similar fashion.
\end{proof}

The first and second inequalities in (\ref{xi-b}) are the so-called global efficiency and reliability bounds, respectively. The reliability bound is used for terminating the adaptive procedure. 
For the ANE method, the global efficiency bound is sufficient for determining the number of new neurons to be added. This is different from the adaptive finite element method in which a local efficiency bound is preferred.

Another estimator of the recovery type is the least-squares (or dual) estimator defined as follows:
 \begin{equation}\label{dual-estimator}
     \eta = \left(\sum_{K\in \cK_n}\eta_K^2\right)^{1/2}
 =\left(\|A^{-1/2}(\hat{\bsigma}_{_\cT}+A\nabla u_{_\cT})\|^2_{0, \Omega} + \|c^{-1/2}(\nabla\!\cdot\! \hat{\bsigma}_{_\cT} +c u_{_\cT} -f)\|^2_{0, \Omega}\right)^{1/2},
 \end{equation}
where $\eta_K$ is the local indicator given by 
 \begin{equation}\label{dual-indicator}
     \eta_K=\left(\|A^{-1/2}(\hat{\bsigma}_{_\cT}+A\nabla u_{_\cT})\|^2_{0, K} + \|c^{-1/2}(\nabla\!\cdot\! \hat{\bsigma}_{_\cT} +c u_{_\cT} -f)\|^2_{0, K}\right)^{1/2}
 \end{equation}
for all $K\in \cK_n$. Let $u$ and $\bsigma = -A\nabla u$ be the solutions of (\ref{primal}) and (\ref{dual}), respectively. Denote the errors by
 \[
 e=u-u_{_\cT} \quad\mbox{and}\quad 
 \be =\bsigma-\hat{\bsigma}_{_\cT}.
 \]
It is easy to see that 
 \[
 \eta^2 =\|A^{-1/2}(\be+A\nabla e)\|^2_{0, \Omega} + \|c^{-1/2}(\nabla\!\cdot\! \be +c\, e)\|^2_{0, \Omega},
 \]
which, together with integration by parts, yields
 \begin{equation}\label{dual-eta}
     |\!|\!|e|\!|\!|^2 + |\!|\!|\be|\!|\!|^2_*
  = \eta^2 + 2\int_{\partial\Omega} e\,(\be\!\cdot\!\bn)\,ds,
 \end{equation}
where $|\!|\!|\cdot|\!|\!|$ and $|\!|\!|\cdot|\!|\!|_*$ are norms given by
 \[
 |\!|\!|v|\!|\!|^2 = \|A^{1/2}\nabla v\|^2_{0, \Omega}
 + \|c^{1/2}v\|^2_{0, \Omega} \quad\mbox{and}\quad |\!|\!|\btau|\!|\!|^2_* =\|A^{-1/2}\btau\|^2_{0, \Omega} + \|c^{-1/2}\nabla\!\cdot\! \btau\|^2_{0, \Omega}.
 \]

\begin{theorem}
Let $\hat{\eta}^2 = \eta^2 + \|u_{_\cT}-g_D\|^2_{1/2,\Gamma_D} +\|\hat{\bsigma}_{_\cT}\cdot\bn+g_N\|_{-1/2,\Gamma_N}^2$. There exists a constant $C_r>0$ 
such that
 \begin{equation}\label{rel-eta-hat}
     |\!|\!|e|\!|\!|\leq C_r\, \hat{\eta}  \quad\mbox{and}\quad
 |\!|\!|\be|\!|\!|_* \leq C_r\, \hat{\eta}.
 \end{equation}
Moreover, if $\|u_{_\cT}-g_D\|_{1/2,\Gamma_D}$ and $\|\hat{\bsigma}_{_\cT}\cdot\bn+g_N\|_{-1/2,\Gamma_N}$ are higher order comparing to $|\!|\!|e|\!|\!|$ and $|\!|\!|\be|\!|\!|_*$, then we have
 \begin{equation}\label{rel-b-1}
     |\!|\!|e|\!|\!|\leq \eta +\mbox{h.o.t} \quad\mbox{and}\quad
 |\!|\!|\be|\!|\!|_* \leq \eta +\mbox{h.o.t}.
 \end{equation}
\end{theorem}

\begin{proof}
 By the definition of the negative norm and the trace inequality, we have 
  \[
  \left|\int_{\Gamma_D} e\,(\be\!\cdot\!\bn)\,ds\right|\leq \|e\|_{1/2,\Gamma_D} \|\be\cdot\bn\|_{-1/2,\Gamma_D} \leq C\, \|e\|_{1/2,\Gamma_D} |\!|\!|\be|\!|\!|_* \leq \dfrac14 |\!|\!|\be|\!|\!|_*^2 + 2C^2 \|e\|_{1/2,\Gamma_D}^2,
  \]
where $C$ is a constant depending on the diffusion and reaction coefficients $A$ and $c$. In a similar fashion, we have 
 \[
 \left|\int_{\Gamma_N} e\,(\be\!\cdot\!\bn)\,ds\right|\leq  \dfrac14 |\!|\!|e|\!|\!|^2 + 2C^2 \|\be\cdot\bn\|_{-1/2,\Gamma_N}^2.
 \]
Hence, we have
 \[
 \left|2\int_{\partial\Omega} e\,(\be\!\cdot\!\bn)\,ds\right| \leq \dfrac12 \big(|\!|\!|e|\!|\!|^2 + |\!|\!|\be|\!|\!|_*^2\big)
 + C\, \left(\|u_{_\cT}-g_D\|^2_{1/2,\Gamma_D} +\|\hat{\bsigma}_{_\cT}\cdot\bn+g_N\|_{-1/2,\Gamma_N}^2
 \right) 
 \]
which, together with (\ref{dual-eta}), implies (\ref{rel-eta-hat}).

If $\|u_{_\cT}-g_D\|_{1/2,\Gamma_D}$ and $\|\hat{\bsigma}_{_\cT}\cdot\bn+g_N\|_{-1/2,\Gamma_N}$ are higher order comparing to $|\!|\!|e|\!|\!|$ and $|\!|\!|\be|\!|\!|_*$, so is $\left|2\int_{\partial\Omega} e\,(\be\!\cdot\!\bn)\,ds\right|$. Now, (\ref{rel-b-1}) is a direct consequence of (\ref{dual-eta}). This completes the proof of the theorem.
\end{proof}

\section{Adaptive Neuron Enhancement (ANE) Method}

Let $u(\bx)$ and $u_{_\cT}(\bx,{\small\btheta}^*_{_\cT})\in \cM^1_n(\tau_k)$ be the solutions of (\ref{primal}) and (\ref{Ritz-n-Q}), respectively. For a given tolerance $\epsilon>0$, this section describes the ANE method (see \cite{LiuCai} for the best least-squares approximation)
to generate a two-layer spline NN such that the approximation accuracy is within the prescribed tolerance, i.e., 
 \begin{equation}\label{tolerance}
     \|u-u_{_\cT}\|_a \leq \epsilon \, \|u_{_\cT}\|_a.
 \end{equation}
For simplicity of presentation, assume that the numerical integration based on a partition $\cT$ is sufficiently accurate (see \cite{LiuCai} on how to adaptively generating a numerical integration mesh).

The procedure of the ANE method is similar to the widely used adaptive mesh refinement (AMR) method for traditional, well-studied mesh-based numerical methods. Unlike the mesh-based methods, the NN method is based on the NN structure determined by the number of neurons in the case of two-layer NNs. This observation suggests that the key question for developing the ANE method is: how many new neurons will be added at each adaptive step? 

To address this question, we introduce the concept of the physical partition of the domain $\Omega$ for a function in $\cM^1_n(\tau_k)$. To this end, let 
\[
u_{_\cT}(\bx,{\small\btheta}^*_{_\cT})=c_0+\sum_{i=1}^n c_i\tau_k(\bomega_i\!\cdot\! \bx -b_i).
\]
The physical partition for $u_{_\cT}(\bx,{\small\btheta}^*_{_\cT})$ is determined by the $n$ hyper-planes $\left\{\bomega_i\!\cdot\! \bx -b_i=0\right\}_{i=1}^{n}$ plus the boundary of the domain $\Omega$, and is denoted by $\cK_n=\{K\}$. Clearly, $\cK_n$ forms a partition of the domain $\Omega$; i.e., the union of all subdomains in ${\cal K}_n$ equals the whole domain $\Omega$ and that any two distinct subdomains of ${\cal K}_n$ have no intersection.


For each element $K\in\cK_n$, denote by $\xi_{_K}$ the local indicator defined in either (\ref{xi}) or (\ref{dual-estimator}). 
We then define a subset $\hat{\cK}_n$ of $\cK_n$ by using either the following average marking strategy:
 \begin{equation}\label{BV-marking-2}
    \hat{\cK}_n =\left\{K\in \cK_n\, :\,
    \xi_{_K}
    \ge \, \dfrac{1}{\#\cK_n}\sum_{K\in \cK_n}\xi_{_K}\right\},
\end{equation}
where $\#\cK_n$ is the number of elements of $\cK_n$, or the bulk marking strategy: finding a minimal subset $\hat{\cK}_n$ of $\cK_n$ such that
\begin{equation}\label{marking-2}
    \sum_{K\in \hat{\cK}_n} \xi^2_{_K}
    \ge \gamma_1\, \sum_{K\in \cK_n} \xi^2_{_K}
    \quad\mbox{for }\,\, \gamma_1\in (0,\,1).
\end{equation}
With the subset $\hat{\cK}_n$, the number of new neurons to be added to the NN is equal to the number of elements in $\hat{\cK}_n$.
With an accurate numerical integration, the ANE method is defined in Algorithm 5.1.

\begin{algorithm}
{\bf {\sc \bf Algorithm 6.1}} Adaptive Neuron Enhancement Method.\\
Given a tolerance $\epsilon >0$ and a numerical integration mesh $\cT$, starting with a two-layer spline NN with a small number of neurons, 
\begin{itemize}
    \item[(1)] solve the optimization problem in (\ref{Ritz-n-Q});
    \item[(2)] estimate the total error by computing $\xi= \left(\sum\limits_{K\in \cK} \xi^2_{_K}\right)^{1/2}$;
    \item[(3)] if $\xi< \epsilon\, \|u_{_\cT}\|_a$, then stop; otherwise, mark $\hat{\cK}_n$ using (\ref{BV-marking-2}) or (\ref{marking-2}), go to Step (4);
    \item[(4)] add $\#\hat{\cK}_n$ neurons to the network, then go to Step (1). 
\end{itemize}
\end{algorithm}

\section{Initialization}

The high dimensional, non-convex optimization problem in (\ref{Ritz-n-Q}) is often solved by iterative optimization methods such as gradient descent (GD), Stochastic GD, Adam, etc. (see, e.g., \cite{BoCuNo2018} for a review paper in 2018 and references therein). 
Usually nonlinear optimizations have many solutions, and the desired one is obtained only if we start from a close enough first approximation. The ANE method provides a natural process for obtaining a good initialization. 

We employ the initialization approach introduced for the best least-squares approximation in \cite{LiuCai}. For readers' convenience, we briefly describe it here.
First, we specify the size of the initial NN and its input and output weights and bias. Input weights and bias are chosen so that the corresponding hyper-planes form a uniform partition of the domain $\Omega$.
The output weights and bias is chosen as the solution of the system of linear equations to be given in (\ref{c^0}).

When the NN is enhanced by adding new neurons in Step (4) of Algorithm (6.1), clearly, parameters corresponding to old neurons will use the current approximation as their initial. Each new neuron is associated with a subdomain $K\in \hat{\cK}_n$ and its initial is chosen so that the corresponding 
hyper-plane passes through the centroid of $K$ and orthogonal to the direction vector with the smallest variance of quadrature points in $K$. This direction vector may be computed by the Principal Component Analysis method (or PCA \cite{PCA1901}).  
For output weights and biases corresponding to new neurons, a simple initial is to set them zero. This means that the initial of the approximation is the current approximation. A better way is to solve problem (\ref{c^0}) for all output weights and bias. 

In the remainder of this section, we describe the system of algebraic equations, that determines the initial of the output weights and bias when the corresponding hyper-planes are fixed. Denote by $\bomega^0=(\bomega_1^0, ... , \bomega_n^0)$ and ${\bf b}^0=(b_1^0, ..., b_n^0)$ the initial of the input weights and bias, respectively. Let  
 \[
u^0_{_\cT}(\bx)=c^0_0 + \sum_{i=1}^n c^0_i \tau_k(\bomega^0_i\cdot\bx - b^0_i)
 \equiv \sum_{i=0}^n c^0_i \varphi_i(\bx;\, \bomega^0_i,\,b^0_i)
\]
be the initial approximation to the solution,  $u_{_\cT}(\bx)\in \cM^1_n(\tau_k)$, of (\ref{Ritz-n-Q}). Then $\bc^0 = (c^0_0, c^0_1, ..., c^0_n)$ is the solution of the following algebraic equations
 \begin{equation}\label{c^0}
     a_{_\cT}(u^0_{_\cT},\varphi_i) =f_{_\cT} (\varphi_i)
     \quad\mbox{for }\,\, i=0,1, ..., n.
 \end{equation}
 
 \begin{lemma}
Assume that the hyper-planes $\{\bomega^0_i\cdot\bx=b^0_i\}_{i=1}^n$ are distinct. 
Then the stiffness matrix ${\bf K} = \Big(a_{_\cT}(\varphi_i(\cdot;\, \bomega^0_i,\,b^0_i),\varphi_j(\cdot;\, \bomega^0_j,\,b^0_j))\Big)_{(n+1)\times (n+1)}$ is symmetric, and positive definite.
\end{lemma}

\begin{proof}
Clearly, ${\bf K}$ is symmetric.
For any $\bv = (v_0,\,v_1,\, ...,\, v_n)^t$,
we have 
\[
\bv^t{\bf K}\bv = a_{_\cT}(v,v),
\]
where $v(\bx) = \sum\limits_{i=0}^n v_i\varphi_i(\bx;\,\bomega^0_i,\, b^0_i)$.
Since $\{\varphi_i\}_{i=0}^n$ are linearly independent (see Lemma~2.1 in \cite{LiuCai}) when the hyper-planes $\{\bomega^0_i\cdot\bx=b^0_i\}_{i=1}^n$ are distinct, $a_{_\cT}(v,v)$ is positive for any nonzero $\bv$, which, in turn, implies that ${\bf K}$ is positive definite.
\end{proof}

\begin{remark}
If there are two hyper-planes are almost linearly dependent, then the stiffness matrix ${\bf K}$ is ill-conditioned even though it is symmetric, positive definite. This is because basis function $\varphi_i(\bx;\,\bomega^0_i,\,b^0_i)=\tau_k(\bomega^0_i\cdot\bx - b^0_i)$ has a non-local support. 
\end{remark}

In one dimension, this difficulty may be overcome by transforming the non-local basis functions to the local nodal basis functions and solving (\ref{c^0}) using the local basis function. More specifically, assume that $b_0^0< b_1^0 < \cdots < b_n^0$. Denote by $h_i=b_{i+1}-b_i$ the length of subinterval $[b_i,b_{i+1}]$, and set 
 \begin{eqnarray*}
  l_i(\bx) &=& h_i^{-1} \varphi_i -\big(h_i^{-1}+h_{i+1}^{-1}\big)\varphi_{i+1} + h_{i+1}^{-1}\varphi_{i+2} \quad\mbox{for }\,\, i=1,..., n-2, \\[2mm]
  l_{n-1}(\bx) &=& h_{n-1}^{-1} \big(\varphi_{n-1}-\varphi_n\big), \quad l_n(\bx) = h^{-1}_n \varphi_n, \quad\mbox{and}\quad l_0(\bx) = \varphi_0-\sum_{i=1}^{n-1} l_i(\bx).
 \end{eqnarray*}
Let $\hat{\bc} =(\hat{c}_0, \hat{c}_1, ..., \hat{c}_n)^t$ be the solution of 
 \[
 \sum\limits_{i=0}^n \hat{c}_j a_{_\cT}(l_j,l_i) =f_{_\cT} (l_i)
     \quad\mbox{for }\,\, i=0,1, ..., n,
 \]
and let $T$ be the linear transformation from $\bvarphi=(\varphi_0, ..., \varphi_n)^t$ to $\bl = (l_0, ..., l_n)^t$. Then 
 \begin{equation} \label{c^01}
 \bc^0 = (c_0^0, ..., c_n^0)^t = T^t \hat{\bc}.    
 \end{equation}



\section{Numerical Experiments}
In this section, we present our numerical results on using the adaptive neuron enhancement (ANE) method to solve diffusion problems based on the Ritz approximation. In all experiments, the minimization problems are iteratively solved by the Adam optimizer \cite{kingma2014adam}. The integrals of the energy functionals are computed numerically using composite mid-point quadrature rules with uniformly distributed quadrature points. For each run during the adaptive enhancement process, the training stops when the value of the energy functional decreases within $0.1\%$ in the last 2000 iterations. 
And the ANE process stops when the user specified accuracy tolerance $\varepsilon$ is obtained, where the error is estimated using the relative recovery error estimator $\xi_{\text{rel}}= \xi /\|\hat{\bsigma}_{_\cT}\|_0 $, where $\|\hat{\bsigma}_{_\cT}\|_0\approx \|u_{_\cT}\|_a$.

\subsection{One-dimensional Poisson Equation}
The first test problem (see \cite{he2018relu,CAI2020,LiuCai}) is a one-dimensional Poisson equation with homogeneous Dirichlet boundary condition defined on the unit interval $\Omega = (0,1)$. For $f(x)=-40000(x^3-2x^2/3+173x/1800+1/300)e^{-100(x-1/3)^2}$, the exact solution of the test problem is given by
\[
 u(x)=x\left(e^{-(x-\frac{1}{3})^2/0.01}-e^{-\frac{4}{9}/0.01}\right). 
 \]

 \begin{figure}[htb!]
  \centering 
   \subfigure[Initial model $u_{_\cT}$ with 10 neurons ${\dfrac{\|u^\prime-u^\prime_{_\cT}\|_0}{\|u^\prime\|_0}}$=0.522380]{
    \label{figpoisson:a} 
    \includegraphics[width=2.5in]{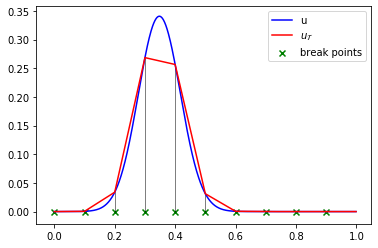}} 
  \hspace{0.3in} 
  \subfigure[Optimized model $u_{_\cT}$ with 10 neurons,  ${\dfrac{\|u^\prime-u^\prime_{_\cT}\|_0}{\|u^\prime\|_0}}$=0.229533]{
    \label{figpoisson:b} 
    \includegraphics[width=2.5in]{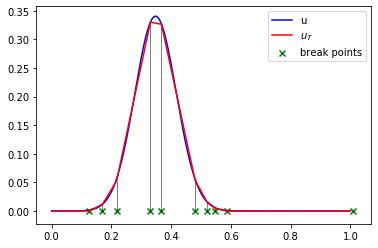}} 
     \hspace{0.3in} 
  \subfigure[Recovered flux $\sigma_{_\cT}$ and the calculated $-u^\prime_{_\cT}$ of 10 neurons,   ${\dfrac{\|\sigma_{_\cT}+u^\prime_{_\cT}\|_0}{\|\sigma_{_\cT}\|_0}}$=0.278647 ]{ 
    \label{figpoisson:c} 
    \includegraphics[width=2.5in]{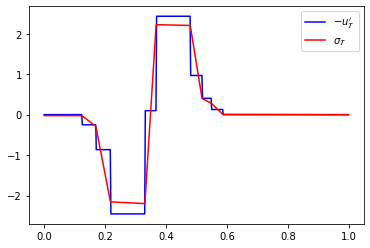}} 
     \hspace{0.3in} 
  \subfigure[Adaptive model $u_{_\cT}$ with 25 neurons, ${\dfrac{\|u^\prime-u^\prime_{_\cT}\|_0}{\|u^\prime\|_0}}$=0.075847]{ 
    \label{figpoisson:d} 
    \includegraphics[width=2.5in]{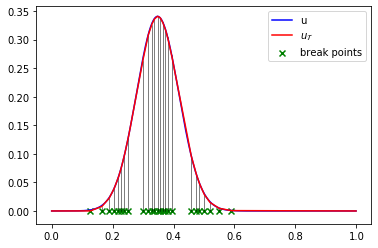}} 
     \hspace{0.3in} 
  \subfigure[Recovered flux $\sigma_{_\cT}$ and the calculated $-u^\prime_{_\cT}$ of 25 neurons, ${\dfrac{\|\sigma_{_\cT}+u^\prime_{_\cT}\|_0}{\|\sigma_{_\cT}\|_0}}$=0.076366]{ 
    \label{figpoisson:e} 
    \includegraphics[width=2.5in]{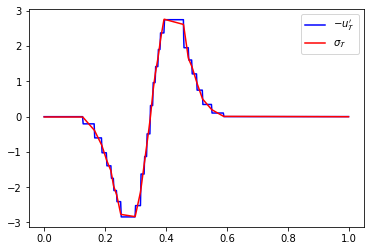}} 
     \hspace{0.3in} 
    \subfigure[A fixed model $u_{_\cT}$ with 25 neurons, ${\dfrac{\|u^\prime-u^\prime_{_\cT}\|_0}{\|u^\prime\|_0}}$=0.151279]{
    \label{figpoisson:f} 
    \includegraphics[width= 2.5in]{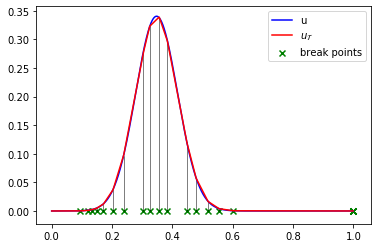}}
     \hspace{0.3in}      
  \caption{Poisson equation approximation results using energy functional as the loss function.} 
  \label{poissonapp} 
\end{figure}


A fixed numerical integration partition $\cal T$ of 1000 uniformly distributed quadrature points is utilized to calculate the energy functional in (\ref{energy}) with $\gamma_D=2000$. For training (optimizing) this primary problem in (\ref{primal}), the learning rate of the Adam optimizer is fixed at 0.002. 

We start from 10 neurons with their breakpoints distributed uniformly and then solve the linear system in (\ref{c^0}) for the initial of the output weights and bias; the initial NN model of $u_{_\cT}$ is depicted in Fig. \ref{figpoisson:a}. 
After optimizing all the parameters in the network, the 10 breakpoints move themselves and the NN outputs an optimized model of non-uniformly distributed breakpoints as shown in Fig. \ref{figpoisson:b}. 

Local error indicator $\xi_K$ is calculated using the recovered $\sigma_{_\cT}$ from $-u^\prime_{_\cT}$ (see Fig. \ref{figpoisson:c} and \ref{figpoisson:e} for a graphical illustration). Elements with large errors are marked by the average marking strategy (see (5.4) in \cite{LiuCai}) and the corresponding neurons are added with proper initialization. This adaptive process repeats itself three runs until our target approximation accuracy $\varepsilon = 0.08$ is reached. Fig.\ref{figpoisson:d} shows the final approximation of adaptive two-layer ReLU NN with 25 neurons.   

For comparison, we also report numerical results using fixed two-layer ReLU NNs with 25 and 50 neurons. Table 1 clearly shows that the accuracy of the adaptive ReLU NN is about the same as that of the fixed NN with twice parameters. 
The approximation of the fixed NN with 25 neurons is depicted in Fig. \ref{figpoisson:f}; the fact that, only 17 out of 25 neurons contribute to the approximation, explains why the fixed NN is not as accurate as the adaptive NN. Finally, we report numerical results from our previous paper \cite{CAI2020} using an over-parametrized DNN of four layers in the last row of Table \ref{poissonnumerical}. Even though the over-parametrized DNN is powerful in approximation, attainable approximation may not be as accurate as that of a proper adaptive/fixed NN with significant less parameters due to the difficulty of non-convex optimization.


\begin{table}[htb]
\caption{Poisson equation: comparing adaptive network with fixed networks using Energy functional }
\label{poissonnumerical}
\centering
\begin{tabular}{  |l |c |c | c | c | c }
	\hline
	NN (hidden layer neurons)
  & $\#$Parameters & ${\dfrac{\|u-{u}_{\tau}\|_0}{\|u\|_0}}$ &  
 ${\dfrac{\|u^\prime-{u^\prime}_{\tau}\|_0}{\|u^\prime\|_0} }$& 
	$\xi_{\text{rel}}=\dfrac{\|\sigma_{_\cT}+u^\prime_{_\cT}\|_0}{\|\sigma_{_\cT}\|_0} $ \\[4mm] \hline
	Fixed 2-layer (25) & 51 &0.012943 &0.149020   & 0.164645 \\ \hline
    Fixed 2-layer (50) & 101 & \textbf{0.006108}  & 0.089470   & 0.095394 \\ \hline	
	Adaptive 2-layer (25) & 51 & 0.007794 & \textbf{0.075847}   & 0.076366 \\ \hline
	Fixed 4-layer (24-14-14) \cite{CAI2020} & 623 & 0.029161 & 0.160666   & - \\ \hline	
\end{tabular}
\end{table}


\subsection{Two-dimensional Poisson Equation with Re-entrant Corner} 

The second test problem is a two-dimensional Poisson equation with pure Dirichlet boundary condition defined on a domain with re-entrant corner $\Omega = \{(r,\theta) |\, r\in (0,1),\,\, \theta \in (0, \frac{3\pi}{2}) \}$. The exact solution 
 \[
 u(r,\theta) = r^{\frac{2}{3}}\sin(\frac{2\theta+\pi}{3}),
 \]
is harmonic, i.e., $\Delta\, u =0$.

 


\begin{table}[htb!]
\caption{Poisson equation with re-entrant corner: comparing adaptive ReLU NN with a fixed NN.}
\label{lshape:num}
\centering
\begin{tabular}{  |l |c |c | c | c | c }
	\hline
	NN (neurons) & $\#$Parameters & ${\dfrac{\|u-{u}_{_\cT}\|_0}{\|u\|_0}}$ &  
 ${\dfrac{\|\nabla (u-{u}_{_\cT})\|_0}{\|\nabla u\|_0} }$& 
	$\xi_{\text{rel}}={\dfrac{\|\bsigma_{_\cT}+\nabla {u}_{_\cT}\|_0 }{\|\bsigma_{_\cT}\|_0}}$ \\[4mm] \hline
    Adaptive 2-layer (20) & 61 & 0.033947 & 0.230226  & 0.233191 \\ \hline
    Adaptive 2-layer (42) & 127 & 0.021939  & 0.154129   & 0.158150  \\ \hline
    Adaptive 2-layer (86) & 259  & \textbf{0.014162} & \textbf{0.064652}  & 0.132546 \\ \hline 
    Fixed 2-layer (86) &259  & 0.025932 & 0.217009  & 0.162507 \\ \hline 

\end{tabular}
\end{table}

 \begin{figure}[htb!]
  \centering
     \subfigure[The exact solution $u$]{ 
    \label{lshape:a} 
    \includegraphics[width=1.9in]{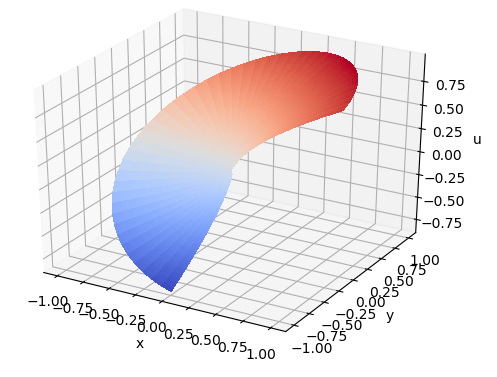}} 
  \subfigure[The exact $\partial_{r} u$]{ 
    \label{lshape:b} 
    \includegraphics[width=1.8in]{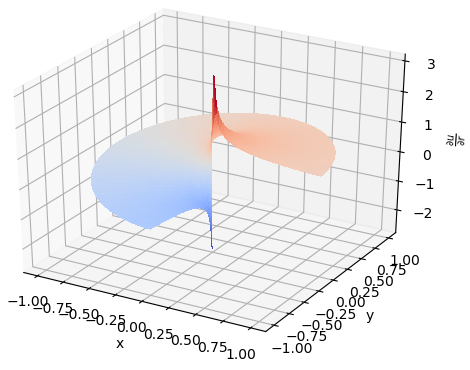}} 
  \subfigure[The exact $\partial_{\theta} u$]{ 
    \label{lshape:c} 
    \includegraphics[width=1.8in]{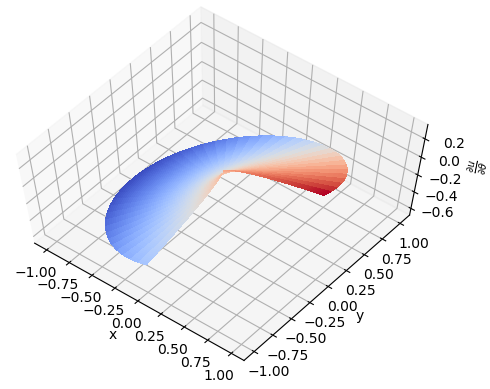}}
   \subfigure[Initial $u_{_\cT}$(20 neurons) ]{ 
    \label{lshape:d} 
    \includegraphics[width=1.8in]{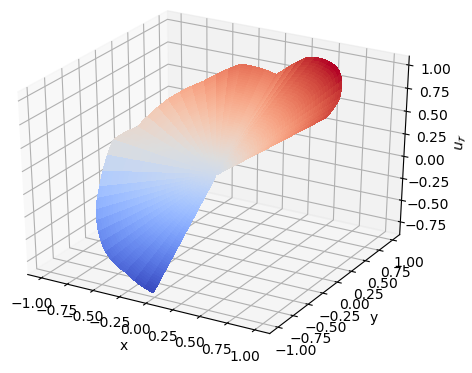}} 
  \subfigure[Initial $\partial_r u_{_\cT}$ (20 neurons)]{ 
    \label{lshape:e} 
    \includegraphics[width=1.8in]{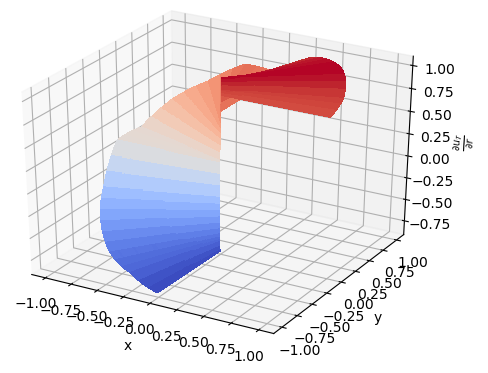}} 
  \subfigure[Initial $\partial_{\theta} u_{_\cT}$ (20 neurons)]{ \label{lshape:f} 
    \includegraphics[width=1.8in]{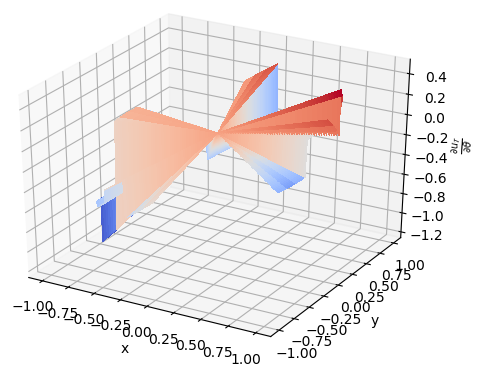}}
   \subfigure[Optimal $u_{_\cT}$(20 neurons)]{ 
    \label{lshape:g} 
    \includegraphics[width=1.8in]{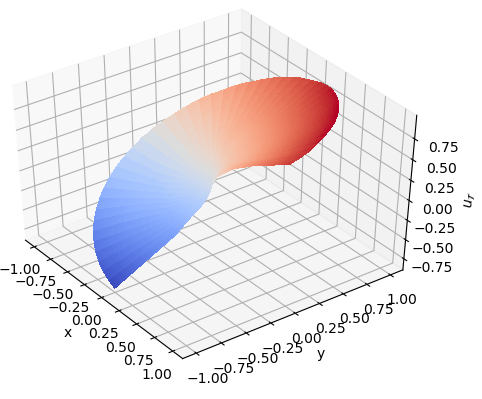}} 
  \subfigure[$\partial_r u_{_\cT}$ (20 neurons)]{ 
    \label{lshape:h} 
    \includegraphics[width=1.8in]{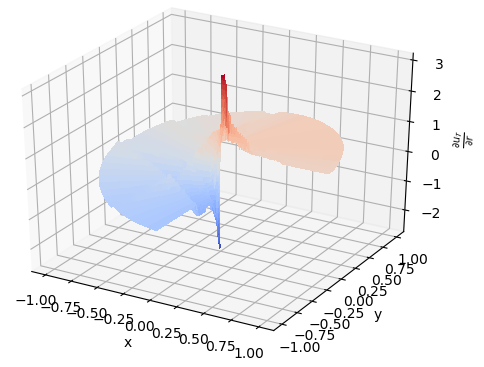}} 
  \subfigure[$\partial_{\theta} u_{_\cT}$ (20 neurons)]{ \label{lshape:i} 
    \includegraphics[width=1.8in]{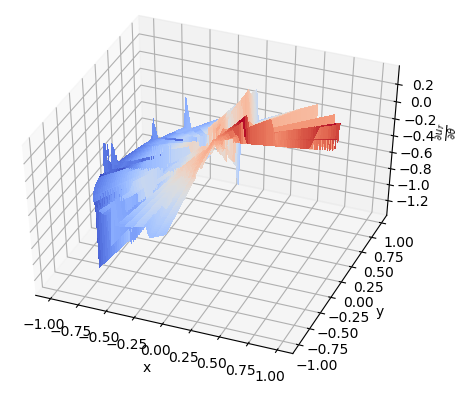}}
    \subfigure[Adaptive NN of $u_{_\cT}$ (86 neurons)]{ 
    \label{lshape:j} 
    \includegraphics[width= 1.8in]{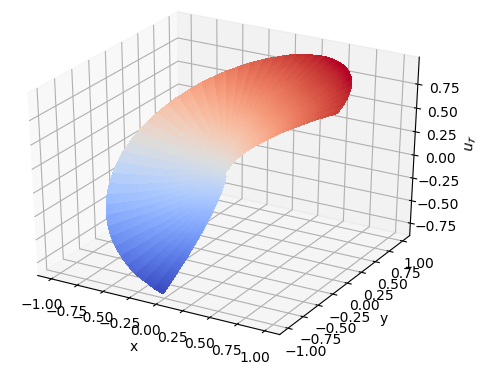}} 
  \subfigure[$\partial_r u_{_\cT}$ (86 neurons) ]{ 
    \label{lshape:k} 
    \includegraphics[width=1.8in]{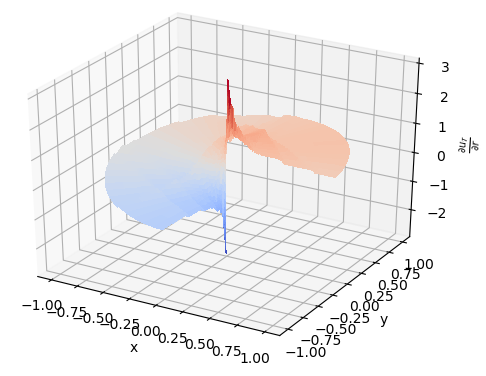}} 
   \subfigure[$\partial_{\theta} u_{_\cT}$ (86 neurons)]{ 
    \label{lshape:l} 
    \includegraphics[width= 1.8in]{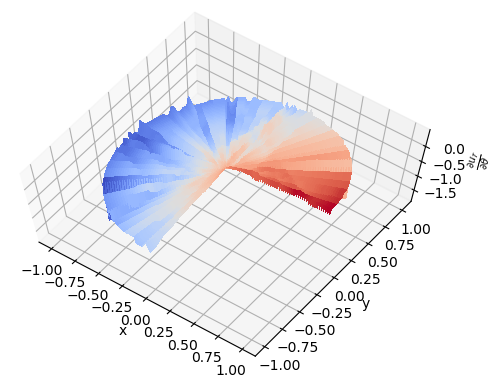}}  
  \caption{Poisson equation with re-entrant corner: Exact solution and results of adaptive two-layer ReLU NN from 20 to 86 neurons.}
  \label{lshape} 
\end{figure}

\begin{figure}[htb]
  \centering
    \subfigure[Initial break lines of 20 neurons]{ 
    \label{lshape2:a} 
    \includegraphics[width=1.75in]{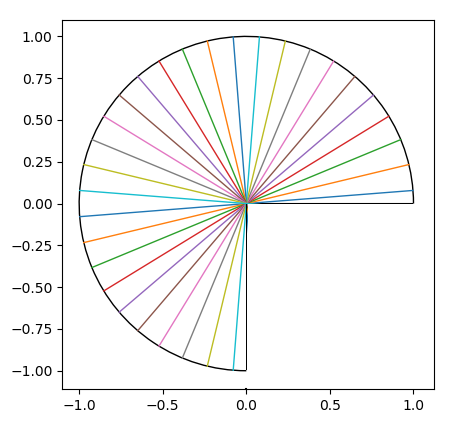}} 
    \hspace{0.1 in}
    \subfigure[Optimal break lines of 20 neurons with marked elements using (\ref{xi})]{ 
    \label{lshape2:b} 
    \includegraphics[width=1.75in]{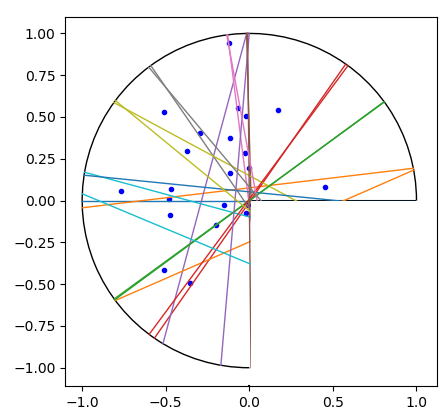}}
    \hspace{0.1 in}
    \subfigure[Elements marked with the exact local error]{ 
    \label{lshape2:c} 
    \includegraphics[width=1.75in]{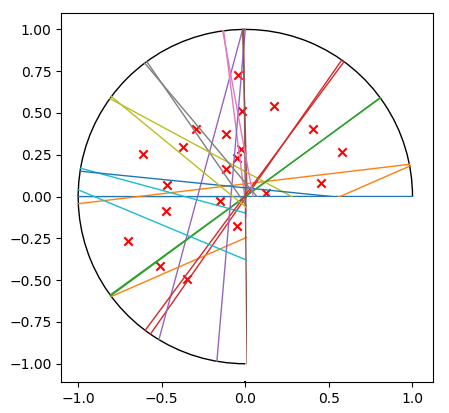}}\\
    \subfigure[Optimal break lines of 42 neurons with marked elements using (\ref{xi})]{ 
    \label{lshape2:d} 
    \includegraphics[width=1.7in]{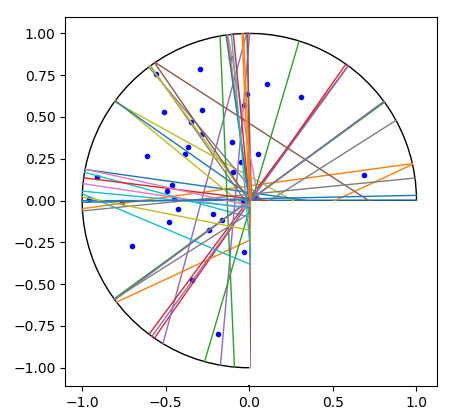}}
    \hspace{0.1in}
    \subfigure[Elements marked with the exact local error]{ 
    \label{lshape2:e} 
    \includegraphics[width=1.75in]{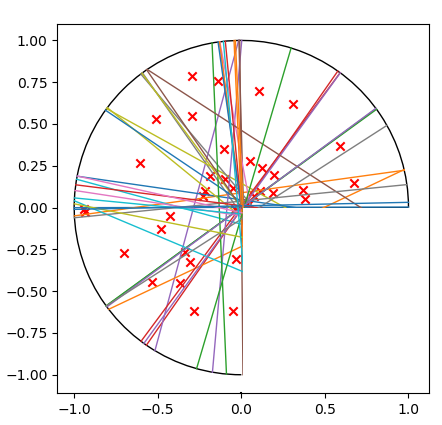}}
    \hspace{0.1in}
    \subfigure[Final break lines of 86 neurons]{ 
    \label{lshape2:f} 
    \includegraphics[width=1.75in]{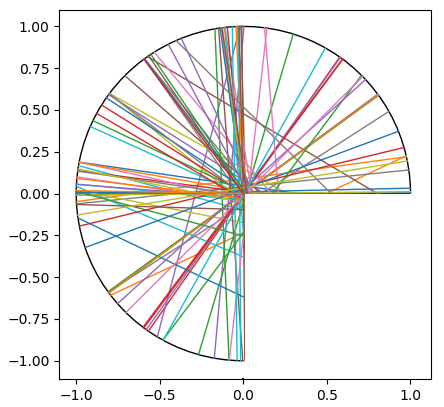}} \caption{Poisson equation with re-entrant corner: break lines generated in the ANE process.}
  \label{lshape2} 
\end{figure}

Numerical integration is calculated using a partition $\cal T$ with $50\times270$ quadrature points uniformly distributed along radial and circumferential directions in a polar coordinate framework. The $\gamma_D$ is set as 200 and the learning rate of the optimizer is fixed at 0.001. The target approximation accuracy is set as $\varepsilon = 0.15$ for this problem considering the difficulty posed by the point singularity at the origin. 

Our adaptive model starts from 20 neurons which are initialized such that the corresponding break lines are distributed uniformly along circumferential direction as shown in Fig. \ref{lshape2:a}. The initial NN model obtained after solving (\ref{c^0}) is illustrated in Fig. \ref{lshape:d}. This initial model gives a fair approximation of $u$ and the relative error in the $L^2$ norm is $0.13$, while approximation to $\nabla u$ (see Fig. \ref{lshape:e} and \ref{lshape:f}) still presents relatively large errors (the relative error in the energy norm is $0.49$). After optimization in the first run, the break lines of these 20 neurons move and form a non-uniform partition of the domain as shown in Fig. \ref{lshape2:b}, which results in an NN model of improved performance, see Table \ref{lshape:num}, first row for the numerical results. The graphical results of $u_{_\cT}$ and $\nabla u_{_\cT}$ approximate by a NN of 20 neurons are depicted in Fig. \ref{lshape:g} - Fig. \ref{lshape:i}. 

During the neuron enhancement step, elements with relative large local error $\xi_K$ are marked using the bulk marking strategy with $\gamma_1 =0.5$ (see (5.5) in \cite{LiuCai}). After two runs of the ANE, adding 22 and 44 neurons respectively, the ANE process stops at 86 neurons with a relative recovery error estimator $\xi_{\text{rel}}= 0.13$. Intermediate results are recorded in Table \ref{lshape:num} and the final visual results of $u_{_\cT}$ and $\nabla u_{_\cT}$ approximated by a NN of 86 neurons are illustrated in Fig. \ref{lshape:j} - Fig. \ref{lshape:l}.

Comparing to a two-layer fixed ReLU NN with the same number of neurons in the hidden layer, the proposed adaptive method converges to a better approximation result (see the last two rows in Table \ref{lshape:num}). This experiment also shows that for a Poisson equation containing a corner singularity, a two-layer ReLU NN is capable of generating a good approximation to the solution. Fig.\ref{lshape2:e} shows the final arrangement generated by the adaptive ReLU NN of 86 neurons in which we observe that the break lines adapt themselves to account for the singular point at the origin. This self-adaptivity of generated physical meshes is a highly desirable feature of using NN model to approximate problems with singularities or discontinuities. 

To verify that the proposed error estimator of the recovery type provides a valid indicator for neuron enhancement, we compare the elements marked using the proposed error indicator (\ref{xi}) and those marked using the exact error, the element marking results for the two intermediate runs in the adaptive process are illustrated in Fig. \ref{lshape2:b}-\ref{lshape2:e}. The comparison results show similar sets of elements being marked for further enhancement, which indicates the validity of the proposed recovery error indicator. 

\subsection{Two Dimensional Interface Problem} 

The third test problem is the intersecting interface problem defined on the unit disk $\Omega = \{(r,\theta) |\, r\in[0,1), \theta \in [0, 2\pi]) \}$ and satisfying 
(\ref{pde}) with $c=f=0$, $\Gamma_N=\emptyset$, and $A = \alpha(\theta)I$, where the diffusion coefficient $\alpha(\theta)$ equals to $1$ in the second and forth quadrants and $R=161.4476387975881$ in the first and third quadrants. 
This is a difficult benchmark test problem for adaptive mesh refinement (see,  e.g., \cite{Morin02, Cai2009}) and the exact solution is $u(r, \theta) = r^{\beta} \mu(\theta)$ with
\[
 \mu(\theta)=
\left\{ \begin{array}{rclll}
\cos((\pi/2-\sigma)\beta)\cdot \cos((\theta-\pi/2+ \rho)\beta), \quad &\text{if} & 0\leq \theta\leq\pi/2,\\[1mm]
\cos(\rho\beta)\cdot \cos((\theta-\pi+ \sigma)\beta), \quad &\text{if} & \pi/2 \leq \theta\leq\pi,\\[1mm]
\cos(\theta\beta)\cdot \cos((\theta-\pi-\rho)\beta), \quad &\text{if} & \pi \leq \theta\leq 3\pi/2,\\[1mm]
\cos((\pi/2-\rho)\beta)\cdot \cos((\theta-3\pi/2-\sigma)\beta), \quad &\text{if} & 3\pi/2 \leq \theta\leq 2\pi.
 \end{array}\right.
\]

 


Considering the inherent difficulty introduced by the intersecting interfaces along x-axis and y-axis and the fact that the recovery estimator over-estimate the true error, we set the stopping criterion as $\xi_{\text{rel}} \leq  \epsilon=0.6$. Numerical integration is calculated on an uniform partition $\cal T$ of $50\times360$ quadrature points. The $\gamma_D$ is set at 200 and a constant learning rate of 0.001 is used throughout the training. For the error estimator of recovery type in (\ref{xi}), we use the identity matrix for $D$. A bulk marking strategy is adopted in the adaptive process with $\gamma_1=0.7$ (see (5.5) of \cite{LiuCai}). 

\begin{figure}[htbp]
  \centering
     \subfigure[Initial break lines (20 neurons)]{ 
    \label{Kellog:a} 
    \includegraphics[width=1.6in]{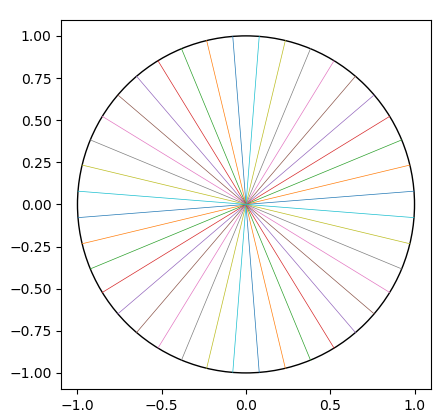}}
  \subfigure[Initial $u_{_\cT}$ (20 neurons)]{ 
    \label{Kellog:b} 
    \includegraphics[width=1.9in]{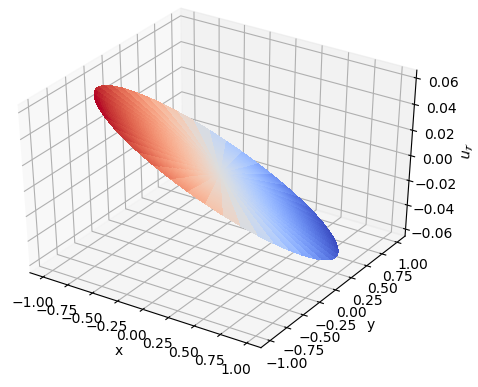}} 
  \subfigure[Initial $-\alpha\partial_xu_{_\cT}$ (20 neurons)]{ 
    \label{Kellog:c} 
    \includegraphics[width=1.8in]{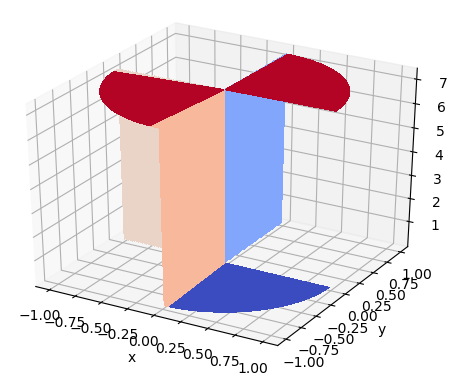}}
    \subfigure[Optimal break lines (20 neurons) with marked elements]{ 
    \label{Kellog:d} 
    \includegraphics[width=1.7in]{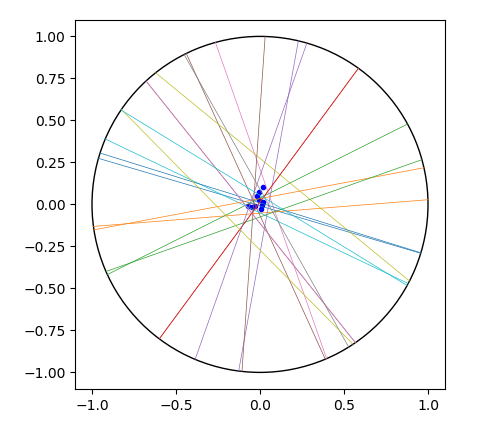}} 
   \subfigure[Optimal $u_{_\cT}$ (20 neurons)]{ 
    \label{Kellog:e} 
    \includegraphics[width=1.7in]{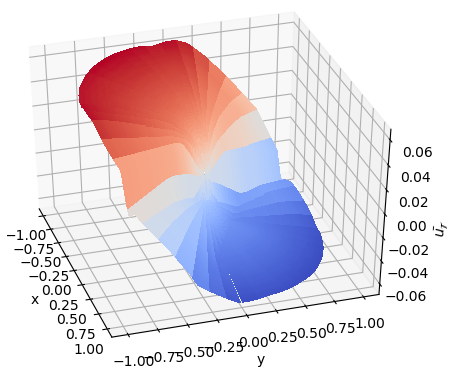}} 
  \subfigure[Optimal $-\alpha\partial_xu_{_\cT}$ (20 neurons)]{ 
    \label{Kellog:f} 
    \includegraphics[width=1.8in]{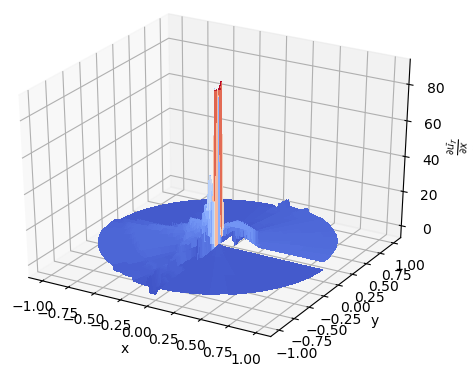}} 
    \subfigure[Final adaptive NN of 150 neurons: break lines]{ 
    \label{Kellog:g} 
    \includegraphics[width= 1.7in]{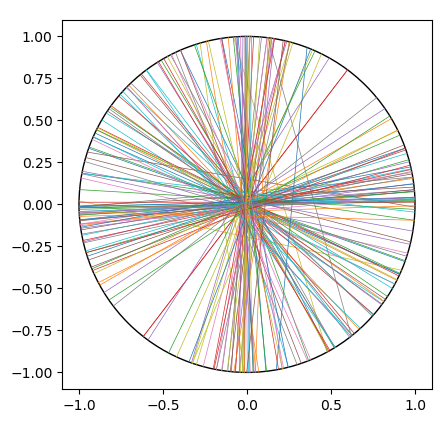}} 
  \subfigure[Adaptive model of $u_{_\cT}$ (150 neurons)]{ 
    \label{Kellog:h} 
    \includegraphics[width=1.8in]{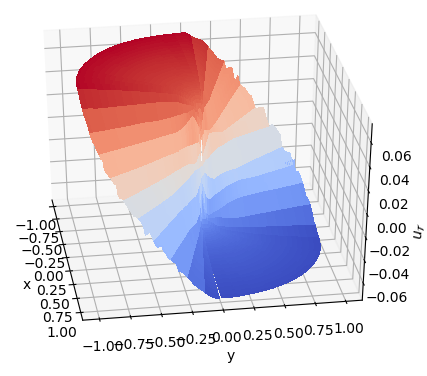}} 
    \subfigure[Final $-\alpha\partial_xu_{_\cT}$ (150 neurons)]{ 
    \label{Kellog:i} 
    \includegraphics[width= 1.9in]{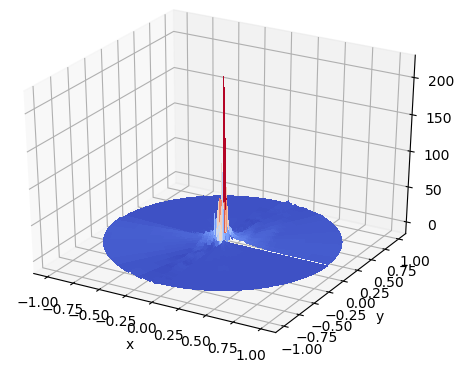}}
      \subfigure[Final $-\alpha\partial_yu_{_\cT}$ (150 neurons) ]{
    \label{Kellog:j} 
    \includegraphics[width=1.8in]{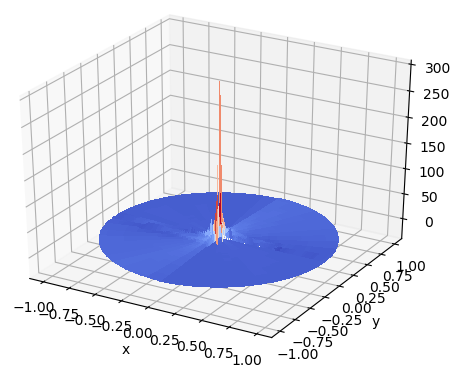}} 
    \subfigure[Recovered flux $\bsigma_{x\cT}$ (150 neurons) ]{ 
    \label{Kellog:k} 
    \includegraphics[width=1.8in]{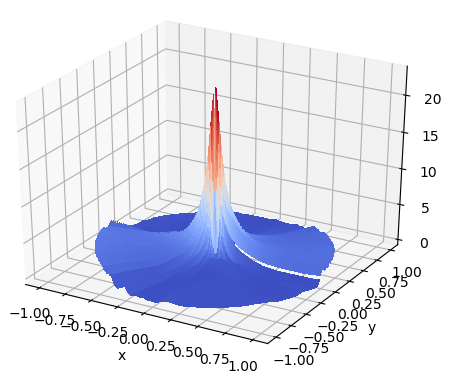}}
    \subfigure[Fixed NN model of $\bar{u}_{_\cT}$ with two hidden layers(20 neurons in each)]{ 
    \label{Kellog:l} 
    \includegraphics[width=1.8in]{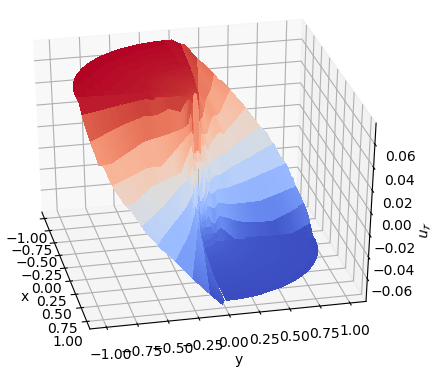}}
  \caption{Kellogg problem: Results of an adaptive 2-layer ReLU NN and a fixed 3-layer ReLU NN.}
  \label{Kellog} 
\end{figure}

Again, we start from a small size NN of 20 neurons, the adaptive process enhances four runs and stops at 150 neurons with the relative error estimator $\xi_{\text{rel}}=0.55 < \epsilon$. The initial NN model, the optimized NN model of 20 neurons and the final model of 150 neurons are all illustrated in Fig. \ref{Kellog}, and the values of the relative error estimator in each run, from 20 neurons to 150 neurons, are recorded in Table \ref{kellog:num}. 

As shown in Table \ref{kellog:num}, the adaptive model of 150 neurons yields a better approximation than the fixed model of the same size.  Comparing to the adaptive finite element method adopted in \cite{Cai2009} which uses more than three thousands of grid points (parameters) in an adaptive refined mesh to achieve a similar result in the relative energy norm, the adaptive ReLU NN presents a more efficient model since all break lines are allowed to move and to adapt to the characteristics of the target function. In this problem, the singularity at the origin and the intersecting interface along axes are captured well through the moving break lines, see Fig. \ref{Kellog:g} for the optimized arrangements generated by the break lines of those 150 neurons.

We notice that the final adaptive model of $u_{_\cT}$ contains certain degrees of oscillation on the boundary (see Fig. \ref{Kellog:h}). To seek a remedy, we further test a fixed three-layer ReLU NN with 20 neurons in each hidden layer. The oscillation is reduced (see Fig. \ref{Kellog:l}) and the relative error in the $L^2$ norm is also improved. The reduction in oscillation is presumably due to the facts that the neurons in deeper layers provide fine scale approximation while those in the first hidden layer provide only coarse scale approximation. Nevertheless, it is also noticed that the adaptive two-layer NN achieves better accuracy in the energy norm comparing to the fixed three-layer NN of similar size (see the last three rows in Table \ref{kellog:num}). This result is perhaps caused by
training of difficult nonlinear optimization problem, because a good initialization is systematically provided through our ANE method for adaptive two-layer NN, but is not available for the fixed three-layer NN. To extend our ANE method to multi-layer NN is our current research project and requires a deeper understanding on the role of depth in a neural network.

\begin{table}[htb]
\caption{Interface Problem: comparing adaptive network with fixed networks using the Ritz formulation.}
\label{kellog:num}
\centering
\begin{tabular}{  |l |c |c | c | c | c }
	\hline
	NN (neurons)
  & $\#$Para. & ${\dfrac{\|u-\bar{u}_{\tau}\|_0}{\|u\|_0}}$ &  
 ${\dfrac{\|A^{1/2}\nabla (u- u_{\tau})\|_0}{\|A^{1/2}\nabla u\|_0} }$& 
	${\xi_{\text{rel}}=\dfrac{\|(\bsigma_{_\cT}+A\nabla u_{_\cT})\|_0 }{\|\bsigma_{_\cT}\|_0}}$ \\[4mm] \hline
    Adaptive (20) & 61 & 0.212283 & 0.985976   & 0.817516 \\ \hline
    Adaptive (31) & 94 & 0.168116 & 0.760418   & 0.749148  \\ \hline  
    Adaptive (52) & 157 & 0.129946 & 0.546362   &0.714189 \\  \hline
     Adaptive (92) & 277 & 0.107181 & 0.362692   &0.634297 \\  \hline   
    Adaptive (150) & 451 &  0.087608 & \textbf{0.047335 }  & 0.549564\\ \hline
    Fixed (150) & 451 & 0.160656  & 0.826836  & 0.717022 \\ \hline  
    Fixed (20-20) & 481  &  \textbf{0.070198} & 0.624581   & 0.535260 \\ \hline     
	
\end{tabular}
\end{table}

\section{Discussion and Conclusion}

To adaptively construct a two-layer spline NN with a nearly minimum number of neurons and parameters such that its approximation accuracy is within the prescribed tolerance, we develop and test the adaptive neuron enhancement (ANE) method for the Ritz approximation to elliptic PDEs in this paper. A key component of the ANE method for its application in PDEs is the development of computable local indicators since the solution of a PDE is unknown. The recovery and the least-squares estimators are introduced. Numerical results for the ReLU NN approximation to problems with corner or intersecting interface singularities show that the recovery estimator is effective. When using other activation functions, the recovery estimator need to be modified by adding weighted $L^2$ norm of the residual of the original equation (see the hybrid estimator in \cite{CaiCai2018, cai2010}) since the recovery estimator may not be reliable. 
The least-squares estimator provides a constant free, guaranteed upper bound of the true error in the energy norm, and hence it is useful to serve as a stopping criterion. 

When a PDE has an underlying minimization principle, our experience suggests that the Ritz formulation is better than various manufactured least-squares formulations as stated in the introduction due to the number of independent variables and the smoothness of the solution. Moreover, a loss functional with fewer independent variables is easier to train than one having more variables. 
Unlike existing NN methods, we approximate the integral of the loss functional by numerical integration. Effect of numerical integration is analyzed for both approximations to a given function and PDE (see Theorem 4.1 of \cite{LiuCai} and Theorem 4.3 of this paper). 

Universal approximation theorem shows that a two-layer ReLU NN is able to accurately approximate any continuous function defined on a compact set in $\R^d$ provided that there are enough neurons in the NN. Indeed, our numerical results demonstrate that problems with corner or intersecting interface singularities may be approximated accurately by the ANE method with fewer degrees of freedom than that of adaptive finite element method. On the other hand, numerical results for problems with sharp circular transition layer \cite{LiuCai} and discontinuous solution \cite{CaiChenLiu2020} show that a three-layer NN is needed in order to approximate the target function well without undesired oscillation. Extension of our ANE method to multi-layer ReLU NN will be presented in a forthcoming article.

\bigskip
\bibliographystyle{elsarticle-num}
\bibliography{Reference}

\begin{thebibliography}{10}
\expandafter\ifx\csname url\endcsname\relax
  \def\url#1{\texttt{#1}}\fi
\expandafter\ifx\csname urlprefix\endcsname\relax\def\urlprefix{URL }\fi
\expandafter\ifx\csname href\endcsname\relax
  \def\href#1#2{#2} \def\path#1{#1}\fi

\bibitem{Berg18}
J.~Berg, K.~Nystrom, A unified deep artificial neural network approach to
  partial differential equations in complex geometries, Neurocomputing 317
  (2018) 28--41.

\bibitem{CAI2020}
Z.~Cai, J.~Chen, M.~Liu, X.~Liu, Deep least-squares methods: An unsupervised
  learning-based numerical method for solving elliptic pdes, Journal of
  Computational Physics 420 (2020) 109707.

\bibitem{Weinan18}
W.~E, B.~Yu, The deep ritz method: A deep learning-based numerical algorithm
  for solving variational problems, Communications in Mathematics and
  Statistics 6~(1) (2018) 1--12.

\bibitem{Karniadakis19}
M.~Raissia, P.~Perdikarisb, G.~Karniadakisa, Physics-informed neural networks:
  A deep learning framework for solving forward and inve, Journal of
  Computational Physics 378 (2019) 686–707.

\bibitem{Sirignano18}
J.~Sirignano, K.~Spiliopoulos, {DGM}: A deep learning algorithm for solving
  partial differential equations, Journal of Computational Physics 375 (2018)
  1139--1364.

\bibitem{LiuCai}
M.~Liu, Z.~Cai, J.~Chen, Adaptive two-layer relu neural network: I. best
  least-squares approximation, submitted (2020).

\bibitem{Xu2020}
J.~Xu, The finite neuron method and convergence analysis, Communications in
  Computational Physics 28 (2020) 1707--1745.

\bibitem{Ciarlet78}
P.~Ciarlet, The finite element method for elliptic problems, Society for
  Industrial and Applied Mathematics, 1978.

\bibitem{Bank}
R.~E. Bank, A.~Weiser, Some a posteriori error estimators for elliptic partial
  differential equations, Mathematics of Computation 44~(170) (1985) 283--301.

\bibitem{cai2010}
Z.~Cai, S.~Zhang, Flux recovery and a posteriori error estimators: conforming
  elements for scalar elliptic equations, SIAM Journal on Numerical Analysis 48
  (2010) 578--602.

\bibitem{Cybenko1989}
G.~Cybenko, Approximation by superpositions of a sigmoidal function,
  Mathematics of Control, Signals, and Systems 2 (1989) 303--314.

\bibitem{HornikS1989}
K.~Hornik, M.~Stinchcombe, H.~White, Multilayer feedforward networks are
  universal approximators, Neural Networks 2 (1989) 359--366.

\bibitem{Petrushev1998}
P.~P. Petrushev, Approximation by ridge functions and neural networks, SIAM
  Journal on Mathematical Analysis 30 (1998) 155--189.

\bibitem{siegel2020approximation2}
J.~W. Siegel, J.~Xu, High-order approximation rates for neural networks with
  $\text{ReLU}^k$ activation functions, arXiv preprint arXiv:2012.07205 (2020).

\bibitem{Weinan20Barron}
W.~E, S.~Wojtowytsch, Some observations on partial differential equations in
  barron and multi-layer spaces (2020).
\newblock \href {http://arxiv.org/abs/arXiv:2012.01484}
  {\path{arXiv:arXiv:2012.01484}}.

\bibitem{Verfurth2013}
R.~Verfurth, A Posteriori Error Estimation Techniques for Finite Element
  Methods, Numerical Mathematics and Scientific Computation, Oxford University
  Press, 2013.

\bibitem{BoCuNo2018}
L.~Bottou, F.~E. Curtis, J.~Nocedal, Optimization methods for large-scale
  machine learning, SIAM Review 60 (2018) 223--311.

\bibitem{PCA1901}
K.~Pearson, On lines and planes of closest fit to systems of points in space,
  The London, Edinburgh, and Dublin Philosophical Magazine and Journal of
  Science 2~(11) (1901) 559--572.

\bibitem{kingma2014adam}
D.~P. Kingma, J.~Ba, Adam: A method for stochastic optimization, arXiv preprint
  arXiv:1412.6980 (2014).

\bibitem{he2018relu}
J.~He, L.~Li, J.~Xu, C.~Zheng, Relu deep neural networks and linear finite
  elements, Journal of Computational Mathematics 38~(3) (2020) 502--527.

\bibitem{Morin02}
P.~Morin, R.~H. Nochetto, K.~G. Siebert, Convergence of adaptive finite element
  methods, Siam Review 44~(4) (2002) 631--658.

\bibitem{Cai2009}
Z.~Cai, S.~Zhang, Recovery-based error estimator for interface problems:
  Conforming linear elements, SIAM Journal on Numerical Analysis 47 (2009)
  2132--2156.

\bibitem{CaiCai2018}
D.~Cai, Z.~Cai, A hybrid a posteriori error estimator for conforming finite
  element approximations, Computer Methods in Applied Mechanics and Engineering
  339 (2018) 320--340.

\bibitem{CaiChenLiu2020}
Z.~Cai, J.~Chen, M.~Liu, Least-squares relu neural network method for linear
  advection-reaction equation, submitted (2020).

\end{thebibliography}

\end{document}